\documentclass{amsart}
\usepackage{amsmath,amssymb}
\usepackage{graphicx}
\usepackage{float}
\usepackage{xcolor}
\usepackage{caption}
\makeatletter
\renewcommand{\fnum@figure}{\figurename~\thefigure.}
\makeatother

\newtheorem{remark}{Remark}

\begin{document}

\title[NS eqs. and diffusion models]{Vortex Stretching in the Navier-Stokes Equations and Information Dissipation in Diffusion Models:
A Reformulation from a Partial Differential Equation Viewpoint}

\author{Tsuyoshi Yoneda} 
\address{Graduate School of Economics, Hitotsubashi University, 2-1 Naka, Kunitachi, Tokyo 186-8601, Japan} 
\email{t.yoneda@r.hit-u.ac.jp} 

\begin{abstract}
    We present a new inverse-time formulation of vortex stretching in the Navier–Stokes equations, based on a PDE framework inspired by score-based diffusion models. By absorbing the ill-posed backward Laplacian arising from time reversal into a drift term expressed through a score function, the inverse-time dynamics are formulated in a Lagrangian manner. Using a discrete Lagrangian flow of an axisymmetric vortex-stretching field, the score function is learned with a neural network and employed to construct backward-time particle trajectories. Numerical results demonstrate that information about initial positions is rapidly lost in the compressive direction, whereas it is relatively well preserved in the stretching direction.
\end{abstract}


\date{\today} 


\maketitle

\section{Introduction: Vortex Dynamics of Navier–Stokes Flows Based on Diffusion Models}

The generation and stretching of vortices in three-dimensional turbulence constitute one of the most fundamental processes in fluid dynamics.
According to classical vortex stretching theory, vorticity is amplified through its interaction with the local strain field, while simultaneously spreading and being smoothed by viscous diffusion.
Such forward-time vortex dynamics have been relatively well understood.
In contrast, the \emph{inverse problem}---namely, how coherent vortex tube structures emerge from a given initial vorticity distribution---remains far from fully understood.
The primary reason lies in the intrinsic irreversibility of the Navier--Stokes equations.

In recent years, score-based diffusion models have rapidly advanced in the field of machine learning, where inverse problems corresponding to forward-time dynamics are successfully formulated.
A particularly noteworthy fact is that the forward-time Fokker--Planck equation governing diffusion models is mathematically identical, up to differences in coefficients, to the vorticity equation describing Burgers vortex tubes.
This structural equivalence provides a theoretical bridge that directly links vortex dynamics with generative processes in diffusion models.
In this study, we exploit this correspondence and construct a new theoretical framework, from the perspective of diffusion models, to understand vortex formation processes in Navier--Stokes flows.
In particular, 
the strain field governs the extent to which the initial vorticity distribution can be reconstructed from the coherent vortex structures observed at later times, that is, the \emph{degree of reversibility of the inverse problem}.
From this viewpoint, the present study raises the following two fundamental questions:
\begin{enumerate}
  \item When vortex dynamics are viewed backward in time, what types of strain fields preserve---or conversely, dissipate---the greatest amount of information about the initial vorticity distribution?
  \item How can the ``information dissipation rate'' associated with a given strain field be quantitatively characterized, and what does it imply about the organization of vortex structures in turbulence?
\end{enumerate}

\noindent
As a first step toward addressing these questions, the present study seeks to gain insight into

\vspace{0.5cm}

\fbox{
  \begin{minipage}{0.9\linewidth}
  how a vortex stretched by a Burgers-type strain field in an axisymmetric Navier--Stokes flow can be pulled back in time, and to what extent the original vorticity distribution can be reconstructed.
  \end{minipage}
}

\vspace{0.5cm}

\noindent
The results provide \emph{a new information-theoretic perspective on vortex dynamics and extend classical vortex stretching theory into a reversible framework inspired by modern diffusion models}.

\section{A Reformulation of Diffusion Models Based on Partial Differential Equation Theory}
In this section, we reformulate diffusion models from a partial differential equation perspective, while keeping probabilistic concepts to a minimum.

\subsection{Axisymmetric vorticity equation with viscosity}

Let us consider a given steady axisymmetric velocity field
$U_r = U_r(r,z)$ and $U_z = U_z(r,z)$,
which satisfies the following divergence-free condition
(the time-dependent case can also be considered, but is omitted here):
\begin{equation*}
\frac{1}{r}\partial_r (r U_r)
+ \partial_z U_z
= 0.
\end{equation*}
Then, the $z$-component of the axisymmetric vorticity equation with viscosity, denoted by $\omega_z$, can be written as
\begin{equation}\label{eq:vorticity}
\begin{split}
\partial_t \omega_z
+ U_r\,\partial_r \omega_z
+ U_z\,\partial_z \omega_z
&=
\omega_z\,\partial_z U_z
+ \omega_r\,\partial_r U_z
+ \nu\left(
\partial_r^2
+ \frac{1}{r}\partial_r
+ \partial_z^2
\right)\omega_z.
\end{split}
\end{equation}
We rewrite the above vorticity equation in the form of a two-dimensional Fokker--Planck equation with a reaction term and an external forcing term:
\begin{equation}
\label{eq:fp-time}
\frac{\partial \rho}{\partial t}
=
-\frac{\partial}{\partial r}\bigl(B_r(r,z)\,\rho\bigr)
-\frac{\partial}{\partial z}\bigl(B_z(r,z)\,\rho\bigr)
+
\nu\left(
\frac{\partial^2 \rho}{\partial r^2}
+
\frac{\partial^2 \rho}{\partial z^2}
\right)
+ S(r,z)\rho
+ F(r,z).
\end{equation}
By comparing Eqs.~\eqref{eq:vorticity} and \eqref{eq:fp-time}, we identify
\begin{equation*}
\label{eq:effective-drift-time}
\rho=\omega_z,\quad
B_r = U_r - \frac{\nu}{r},
\quad
B_z = U_z,\quad
S = \partial_r B_r + 2\,\partial_z B_z + \frac{\nu}{r^2},
\quad
F = \omega_r\,\partial_r U_z.
\end{equation*}
In the case of a Burgers-type strain field, we have
$U_r = -a r$ and $U_z = 2 a z$ $(a>0)$, which yields $F = 0$.

\subsection{Derivation of the time-reversed PDE}
Here we set $f=(B_r,B_z)$ and $x=(x_1,x_2)=(r,z)$.
To match the notation used in the implementation, we define $\sigma:=\sqrt{2\nu}$. 
For simplicity, we consider the case without external force, i.e.\ $F=0$.
Then the density $\rho(t,x)$ satisfies
\begin{equation}
  \partial_t \rho(t,x)
  = -\nabla_x\cdot\bigl(f(x)\,\rho(t,x)\bigr)
    + \frac{\sigma^2}{2}\,\Delta_x \rho(t,x)+S(x)\rho(t,x),
  \qquad t>0.
  \label{eq:FP-forward}
\end{equation}
This is the forward Fokker--Planck equation.
If one naively reverses time in this equation and formulates the corresponding
\emph{inverse problem} by straightforward calculations, the resulting PDE for the backward ``time''
contains a backward Laplacian explicitly and is therefore clearly ill-posed.
The key to properly formulating the inverse problem is to rewrite this backward Laplacian
as part of a drift term.
First, fix a terminal time $T>0$ and define
\begin{equation*}
  \tau = T - t,
  \qquad
  \rho^*(\tau,x) := \rho(T-\tau,x).
\end{equation*}
Then
\begin{equation*}
  \partial_\tau \rho^*(\tau,x)
  = -\partial_t \rho(t,x)\big|_{t=T-\tau}.
\end{equation*}
Substituting \eqref{eq:FP-forward}, we obtain
\begin{align}
\nonumber
  \partial_\tau \rho^*(\tau,x)
  &= -\Bigl(
       -\nabla_x\cdot\bigl(f(x)\,\rho(t,x)\bigr)
       + \frac{\sigma^2}{2}\,\Delta_x \rho(t,x)
       +S(x)\rho(t,x)
     \Bigr)_{t=T-\tau}
\\
  &= \nabla_x\cdot\bigl(f(x)\,\rho^*(\tau,x)\bigr)
     - \frac{\sigma^2}{2}\,\Delta_x \rho^*(\tau,x)+S(x)\rho^*(\tau,x),
  \label{eq:FP-time-reversed}
\end{align}
which is the time-reversed form of the equation.
On the other hand, a non-diffusive equation with drift $\tilde{b}(\tau,x)$
in the backward time variable $\tau$ can be written as
\begin{equation}
  \partial_\tau \rho^*(\tau,x)
  = -\nabla_x\cdot\bigl(\tilde{b}(\tau,x)\,\rho^*(\tau,x)\bigr)
  +S(x)\rho^*(\tau,x).
  \label{eq:FP-backward-candidate}
\end{equation}
Since \eqref{eq:FP-time-reversed} and \eqref{eq:FP-backward-candidate}
are both PDEs for the same $\rho^*(\tau,x)$, their right-hand sides must coincide.
Thus we require
\begin{equation}\label{eq:div-equality}
  \nabla_x\cdot\bigl(f\,\rho^*\bigr)
  - \frac{\sigma^2}{2}\,\Delta_x \rho^*
  \;=\;
  -\nabla_x\cdot\bigl(\tilde{b}\,\rho^*\bigr).
\end{equation}
We would like to express the left-hand side in divergence form.
A standard calculation shows that by setting
\begin{equation*}
  \tilde{b}\rho^*
  = -f\rho^* + \frac{\sigma^2}{2}\nabla_x \rho^*,
\end{equation*}
equation \eqref{eq:div-equality} is satisfied.
Hence, in regions where $\rho^*(\tau,x)>0$, the drift can be written as
\begin{equation}  \label{eq:drift-reversed}
  \tilde{b}(\rho_0,t,x)
  = -f(x) + \frac{\sigma^2}{2}\,\frac{\nabla_x \rho^*(t,x)}{\rho^*(t,x)}
  = -f(x) + \frac{\sigma^2}{2}\,\nabla_x \log \rho^*(t,x).
\end{equation}
This drift term arises from the backward Laplacian and explicitly depends
on the initial density $\rho_0$.
The quantity $\nabla_x \log \rho^*(t,x)$ is referred to as the
(time-reversed) \emph{score function}.
In terms of PDEs, we may write
\begin{equation*}
  \partial_\tau\rho^*
  = -\nabla\cdot\left(
  \tilde b\rho^*\right)
  +S\rho^*.
\end{equation*}
Collecting these results, we see that the problem ultimately reduces to solving the following system of PDEs, expressed in terms of the Lagrangian flow:
\begin{eqnarray}
\label{velocity}
  \partial_\tau(\rho^*\circ\Phi)
  &
  =& 
  \left(-\nabla\cdot \tilde b
+S\right)(\rho^*\circ\Phi),\\
\label{lagrangian}
\partial_\tau\Phi
&
=&
\tilde b
\circ\Phi.
\end{eqnarray}
Recall that the drift $\tilde b$ depends on $\rho_0$, and therefore forms a family
of functions parameterized by the initial condition.
In principle, one must address the question of whether there exists an initial
density $\rho_0$ corresponding to a given terminal density $\rho(T,\cdot)$.
This is a difficult and highly nontrivial problem.
Rather than tackling this problem directly, we adopt machine learning as a
breakthrough approach: for each $(x,t)$, an appropriate value is selected
from the family of initial-condition-dependent functions by learning.

\begin{remark}
As emphasized repeatedly, a purely deterministic manipulation of the equations yields
\begin{eqnarray*}
  \partial_\tau(\rho^*\circ\Phi)
  &
  =& 
  \left(-\frac{\sigma^2}{2}\Delta+\nabla\cdot f
+S\right)(\rho^*\circ\Phi),\\
\partial_\tau\Phi
&
=&
-f
\circ\Phi,
\end{eqnarray*}
which explicitly involves a backward Laplacian.
Such a PDE is therefore ill-posed and mathematically meaningless.
\end{remark}

\subsection{Reformulation of the score function via a discrete Lagrangian flow representation of the heat equation solution}

To begin with, for simplicity, we examine the score function generated by the heat equation
and its discrete Lagrangian flow representation.
The PDE analysis employed here is entirely standard.
In essence, the diffusion term is represented by the heat kernel,
and the order of the time parameter $t$ is evaluated based on parabolic scaling.
The heat equation on $\mathbb{R}^d$ is given by
\begin{equation*}
  \partial_t \rho(t,x) = \nu \Delta_x \rho(t,x),
  \qquad t>0,\ x\in\mathbb{R}^d,
\end{equation*}
its fundamental solution is
\begin{equation*}
  G(t,x)
  =
  \frac{1}{(4\pi\nu t)^{d/2}}
  \exp\Bigl(
    -\frac{|x|^2}{4\nu t}
  \Bigr),
\end{equation*}
which can be regarded as the solution with a Dirac delta initial condition.
When considering a discrete Lagrangian flow, it suffices to focus on the solution at time
$\Delta t$ corresponding to this Dirac delta initial data.
In other words, we wish to estimate the location $\tilde{x}$ at time $\Delta t$
of a fluid particle that was located at position $x$ at the initial time.
This is described by the following probability density (i.e., the solution $\rho$ of the heat equation):
\begin{equation*}
\rho(\Delta t,\tilde x)=
  \frac{1}{(4\pi\nu\Delta t)^{d/2}}
  \exp\Bigl(
    -\frac{|\tilde{x}-x|^2}{4\nu\Delta t}
  \Bigr).
\end{equation*}
This expression can be naturally rewritten in the form of a discrete Lagrangian flow as
\begin{equation*}
  \tilde x = x + \sqrt{2\nu\Delta t}\,\varepsilon,
  \qquad
  \varepsilon \sim \mathcal{N}(0, I),
\end{equation*}
where $\varepsilon$ denotes Gaussian noise.
Accordingly, in the discrete Lagrangian flow the velocity $\tilde b$ 
is represented as
\begin{equation*}
    b(x,\tilde x):=\frac{x-\tilde x}{\Delta t}
    =2\nu\nabla_{\tilde x}\log\rho(\Delta t,\tilde x)
    =-\sqrt{\frac{2\nu}{\Delta t}}\varepsilon .
\end{equation*}
For a fixed target point $x$, $\tilde b(x,\tilde x)$ represents the backward-in-time Lagrangian velocity that pulls a particle located at $\tilde x$ after time $\Delta t$ toward the most probable position $x$ at the previous time.
The second equality
\begin{equation*}
    2\nu\nabla_{\tilde x}\log\rho(\Delta t,\tilde x)
    =-\sqrt{\frac{2\nu}{\Delta t}}\varepsilon
\end{equation*}
interprets the Gaussian noise appearing in the forward-time particle trajectory and the maximum-likelihood gradient of the backward-time probability density as arising from the same Gaussian transition (the most important bridge in reformulating the diffusion model from a PDE perspective).
From this expression, we further obtain
\begin{equation*}
    \nabla_{\tilde x}\cdot\tilde b(x,\tilde x)
    =-\frac{d}{\Delta t},
\end{equation*}
which diverges to $-\infty$ in the continuous limit $\Delta t\to 0$.
While the backward-time process of the heat equation is ill-posed in the continuous-time limit, temporal discretization makes the reverse dynamics explicit.
 In this sense, temporal discretization is not merely a numerical approximation but an essential operation for visualizing the backward process.
However, if one FORMALLY inserts this $\tilde b$ into a Lagrangian flow by taking $\Delta t\to 0$, one obtains
\begin{equation*}
    \partial_{\tau}\Phi
    =
    \tilde b(x,\cdot)\circ\Phi,
\end{equation*}
which is ill-defined, since the destination point $x$ in the backward-time direction appears explicitly as a variable in the velocity field.
This difficulty motivates the introduction of machine learning methods to overcome the dependence on $x$.

We next extend this sequence of calculations to equations that include
both advection and reaction terms.

\subsection{Discrete Lagrangian flow representation and score function formulation for solutions of the Fokker–Planck equation with reaction terms}

Let $p(t,x,y)$ denote the fundamental solution of the following
Fokker--Planck equation with a reaction term
(for simplicity we consider the one-dimensional case, but the extension
to higher dimensions is straightforward):
\begin{equation}
\label{eq:fp}
\begin{cases}
\displaystyle
\partial_t p(t,x,y)
=
-\partial_y\!\bigl(b(y)\,p(t,x,y)\bigr)
+\frac{1}{2}\,\partial_y^2\!\bigl(a(y)\,p(t,x,y)\bigr)
+S(y)p(t,x,y),
\\[6pt]
p(0,x,y)=\delta(y-x).
\end{cases}
\end{equation}
Here $x$ is a spatial parameter representing the initial position,
and the unknown function is $p$ viewed as a function of $y$.
In this section, we employ a standard analytical technique for studying
the short-time behavior of fundamental solutions of parabolic equations,
namely,
{\bf ``a time-direction Taylor expansion of the semigroup representation
via the weak formulation''}.
This approach justifies the fact that the true fundamental solution
and a Gaussian kernel share the same first-order short-time expansion.
Multiplying \eqref{eq:fp} by a test function $\varphi(y)$
and integrating by parts with respect to $y$, we obtain
\begin{equation}
\label{eq:weak_fp}
\frac{d}{dt}\int p(t,x,y)\varphi(y)\,dy
=
\int_{\mathbb{R}} p(t,x,y)\,
\Bigl(
b(y)\varphi'(y)
+\frac{1}{2} a(y)\varphi''(y)
+S(y)\varphi(y)
\Bigr)\,dy.
\end{equation}
We define the differential operator appearing on the right-hand side by
\[
(L\varphi)(y)
:=
b(y)\varphi'(y)
+\frac{1}{2} a(y)\varphi''(y)
+S(y)\varphi(y).
\]
From the initial condition we have
\[
\int_{\mathbb{R}} \delta(y-x)\,\varphi(y)\,dy
=
\varphi(x),
\]
where the integral is understood in the sense of distributions.
Evaluating \eqref{eq:weak_fp} at $t=0$, we obtain
\[
\begin{aligned}
\frac{d}{dt}\int p(t,x,y)\varphi(y)\,dy\bigg|_{t=0}
&=
\int_{\mathbb{R}} \delta(y-x)\,(L\varphi)(y)\,dy
\\
&=
(L\varphi)(x)
=
b(x)\varphi'(x)
+\frac{1}{2} a(x)\varphi''(x)
+S(x)\varphi(x).
\end{aligned}
\]
Therefore, by taking the Taylor expansion with respect to $t$ at $t=0$,
we obtain
\begin{equation}
\label{eq:weak_expansion}
\int_{\mathbb{R}} p(\Delta t,x,y)\,\varphi(y)\,dy
=
\varphi(x)
+\Delta t
\Bigl(
b(x)\varphi'(x)
+\frac{1}{2} a(x)\varphi''(x)
+S(x)\varphi(x)
\Bigr)
+o(\Delta t).
\end{equation}
Next, we consider the constant-coefficient equation obtained by freezing
the coefficients at the point $x$:
\begin{equation*}
\label{eq:frozen_fp}
\partial_t q
=
-b(x)\,\partial_y q
+\frac12 a(x)\,\partial_y^2 q
+S(x)q,
\qquad
q(0,y)=\delta(y).
\end{equation*}
Since this is a linear equation, standard calculations show that its
fundamental solution is given explicitly by
\[
q(t,y)
=
e^{S(x)t}\frac{1}{\sqrt{2\pi a(x)t}}
\exp\!\left(
-\frac{(y-b(x)t)^2}{2a(x)t}
\right).
\]
Therefore,
\begin{equation}
\label{eq:gaussian_expansion}
\begin{split}
&\int_{\mathbb{R}} q(\Delta t,y-x)\varphi(y)\,dy\\
=&
e^{S(x)\Delta t}
\int_{\mathbb{R}}
\varphi\bigl(x+b(x)\Delta t+\sqrt{a(x)\Delta t}\,z\bigr)
\frac{e^{-z^2/2}}{\sqrt{2\pi}}\,dz\\
=&
\varphi(x)
+\Delta t
\Bigl(
b(x)\varphi'(x)
+\frac12 a(x)\varphi''(x)
+S(x)\varphi(x)
\Bigr)
+O(\Delta t^{3/2}).
\end{split}
\end{equation}
Note that the terms of order $\Delta t^{1/2}$ vanish because
$z e^{-z^2/2}$ is an odd function. This is a key point shared with
It\^o integrals and the present weak approximation.
We now replace the test function $\varphi$ by a sequence of functions
converging to the Dirac delta.
Taking into account the competition between this sequence and the
$O(\Delta t^{3/2})$ term, we choose the sequence so that
$\Delta t^{3/2}$ converges to zero more rapidly.
By comparing \eqref{eq:weak_expansion} and \eqref{eq:gaussian_expansion},
this justifies that the fundamental solution $p(\Delta t,x,\cdot)$
and the Gaussian kernel $G$ share the same first-order short-time expansion.

As in the case of the heat kernel discussed in the previous section,
we obtain the discrete Lagrangian flow
\begin{equation*}
\label{eq:discrete_update}
\tilde x
=
x
+
b(x)\,\Delta t
+
\sqrt{a(x)\,\Delta t}\;\varepsilon,
\qquad
\varepsilon\sim\mathcal N(0,1),
\end{equation*}
where $\varepsilon$ again denotes Gaussian noise.
The reformulation of the score function can be derived in the same manner
as for the heat kernel.
To bring the discussion closer to implementation, we now set
$a=\sigma^2$ to be a constant.
The solution of the PDE at time $\Delta t$ is then given by
\begin{equation*}
\label{eq:local-weight}
p(\Delta t,x,\tilde x)
=
\frac{1}{\sqrt{2\pi\sigma^2\Delta t}}
\exp\bigl(S(x)\Delta t\bigr)
\exp\!\Bigl(
-\frac{
\bigl|
\tilde x - x - b(x)\,\Delta t
\bigr|^2
}{
2\sigma^2\Delta t
}
\Bigr)
+o(\Delta t).
\end{equation*}
From this representation formula, we obtain the following score function:
\begin{align*}
  \nabla_{\tilde x} \log p(\Delta t,x,\tilde x)
  &=
  -\frac{
       \tilde x - x - b(x)\,\Delta t
     }{
       \sigma^2\Delta t
     }.
  \label{eq:local-score}
\end{align*}
Therefore, we can appropriately define the following discrete version of the drift $\tilde b$, which explicitly incorporates the backward-time contribution:
\begin{equation*}
    \tilde b(x,\tilde x)
    := \frac{x - \tilde x}{\Delta t}
    = - b(x) + \sigma^2 \nabla_{\tilde x} \log p(\Delta t, x, \tilde x).
\end{equation*}


\begin{remark}
Summarizing the insights obtained so far, the discrete Lagrangian flow
for the viscous axisymmetric vorticity equation is given by
\begin{equation*}
\label{eq:sde-axi-time}
\boxed{
\begin{aligned}
R_{k+1} &= \left(U_r(R_k,Z_k) - \frac{\nu}{R_k}\right)\,\Delta t
+ \sigma\sqrt{\Delta t}\,\varepsilon^{(r)},\\[4pt]
Z_{k+1} &= U_z(R_k,Z_k)\,\Delta t
+ \sigma\sqrt{\Delta t}\,\varepsilon^{(z)},
\end{aligned}
}
\end{equation*}
where $\varepsilon^{(r)}$ and $\varepsilon^{(z)}$ are independent
Gaussian random variables.
\end{remark}

\subsection{Score functions arising in backward discrete Lagrangian flows and machine learning}
At each discrete step $k$, a simple algebraic manipulation of the score function
derived in the previous sections
yields the following recurrence relation:
\begin{equation*}
  x_k
  =
  x_{k+1}
  +
  \tilde b(x_k,x_{k+1})\,\Delta t.
\end{equation*}
This provides a discrete counterpart of the backward-in-time Lagrangian flow
\eqref{lagrangian}.
However, the objective of this study is to express the backward Lagrangian flow
in a closed form that depends only on $x_{k+1}$.
To this end, we approximate $\tilde b(x_k,x_{k+1})$ by a function of the form
\begin{equation*}
  v(x_{k+1},k).
\end{equation*}
Here, $v$ is a function conditioned on the time step $k$ and is represented
in practice by a neural network.


\begin{remark}
In standard diffusion models, a recurrence relation of the form
\begin{equation*}
x_t = \sqrt{1-\beta_t}\,x_{t-1} + \sqrt{\beta_t}\,\varepsilon_t,
\quad
\beta_t\in(0,1),\quad
\varepsilon_t \sim \mathcal{N}(0, I)
\end{equation*}
is commonly assumed.
Defining
\begin{equation*}
\alpha_t := 1 - \beta_t,\quad
\bar{\alpha}_t := \prod_{s=1}^{t} \alpha_s,
\end{equation*}
a direct calculation (using the independence of the Gaussian noises
$\varepsilon_t$) yields
\begin{equation*}
x_t = \sqrt{\bar\alpha_t}\,x_{0} + \sqrt{1-\bar\alpha_t}\,\varepsilon,
\qquad
\varepsilon \sim \mathcal{N}(0, I).
\end{equation*}
This expression provides an explicit representation of $x_t$ in terms of the
initial position $x_0$ and Gaussian noise $\varepsilon$, and is therefore
highly convenient.
In the present setting, however, the Navier--Stokes equations constitute the
fundamental underlying model, and this representation is no longer applicable.
\end{remark}
Now suppose that $N$ sample trajectories
\[
  \{ x_k^{(n)} \}_{k=0}^{L-1},
  \qquad n=1,\dots,N,
\]
obtained from the forward discrete Lagrangian flow are given.
To evaluate the approximation accuracy of $v$, we introduce the following
empirical functional:
\begin{equation*}
  \mathcal{J}[v]
  :=
  \frac{1}{N}\sum_{n=1}^N
  \sum_{k=0}^{L-1}
  \left|
    v\bigl(x_{k+1}^{(n)},k\bigr)
    -
    \tilde b\bigl(x_k^{(n)},x_{k+1}^{(n)}\bigr)
  \right|^2 .
\end{equation*}
By minimizing this functional, we obtain
\begin{equation*}
  v^\ast
  :=
  \arg\min_v \mathcal{J}[v],
\end{equation*}
and the corresponding backward discrete Lagrangian flow is given by
\begin{equation*}
  \boxed{
  x_k
  =
  x_{k+1}
  +
  v^\ast(x_{k+1},k)\,\Delta t
  }.
\end{equation*}
In this way, the generative dynamics in the backward direction are constructed
by using the score function learned from the forward discrete Lagrangian flow.

\section{Implementation of backward-time dynamics for axisymmetric Navier–Stokes flows}
In this section, we present the learning results.

\subsection{Hyperparameter and evaluation metric settings}

We first describe the various hyperparameter settings.
Here, we implemented the Burgers-type strain field with
$\nu = 1$ and $0.01$, and
$U_r = -a r$, $U_z = 2 a z$ with $a=1$.
The time interval $[0,T]$ was divided into $L$ steps,
with time step size $\Delta t = T/(L-1)$.
In the present implementation, we set $T=2$ and $L=200$.
The initial positions of fluid particles were specified using a scale
parameter $s>0$ as
\begin{equation*}
(R_0, Z_0) = (e^{2aT}s, s),
\end{equation*}
and we examined the cases $s = 1,2,\cdots,11,12$.
This choice ensures that the scale of $R$ at the initial time
matches that of $Z$ at the terminal time.
In generating the forward discrete Lagrangian flow (training data),
trajectories in which $R$ became negative at any intermediate step
were discarded.
Training data were generated so that the total number of samples
was $N=10{,}000$, retaining only trajectories with $R>0$
(8{,}000 samples for training and 2{,}000 for validation).

The initial-value reconstruction error was evaluated using a
relative mean absolute error (MAE),
normalized by the total displacement of each component.
For example, for the $R$ component, we used
\begin{equation}
\frac{1}{N}\sum_{i=1}^N
\frac{|R_0^{(i)} - \hat{R}_0^{(i)}|}
     {\sum_{k=0}^{L-2} |R_{k+1}^{(i)} - R_k^{(i)}|},
\end{equation}
where $\hat R_0^{(i)}$ is the corresponding prediction.
The same definition was applied to the $Z$ component.
For each scale parameter $s$, multiple independent trials were conducted,
and the mean and standard deviation of the relative MAE were computed.
To enhance the reliability of the learning results, trials were repeated
until the relative standard error of the mean MAE
(the standard error normalized by the mean)
was reduced to approximately 6\%, which required about 80 runs.

\subsection{Architecture}

In this study, we adopted the framework of score-based generative models
based on continuous-time stochastic differential equations
proposed by Song \emph{et al.}~\cite{Song2019,Song2021}.
Specifically, we employed an {\bf \emph{orthodox implementation of the score
network and adopted a minimal architecture retaining only the elements
required for the present problem setting}}.
Compared with high-dimensional models possessing complex spatial structures,
such as those used in AI-based image generation,
the system considered here has an extremely low-dimensional structure,
tracking a discrete Lagrangian flow in two dimensions $(R,Z)$.
Therefore, adopting such a minimal configuration is well justified.

The score network architecture consists of a multilayer perceptron
that takes the state variables $(R,Z)$ and time $t$ as inputs.
Temporal information is transformed into a latent representation
using a standard time-embedding technique based on sine and cosine functions.
This time-embedded representation is concatenated with the latent
representation of the state variables $(R,Z)$ and fed into the main network,
which consists of multiple fully connected layers and activation functions.
As a result, a score function representation that is continuous
in both time and space is obtained.
The output corresponds to the score for each component $(R,Z)$,
and a single network simultaneously learns the scores
over all time steps $k = 0,1,\dots,L-1$.

\subsection{Learning results ($\nu=1$ and $\nu=0.01$)}

The learning results are shown below.

\begin{figure}[H]
\centering
\begin{minipage}{0.48\textwidth}
  \centering
  \includegraphics[width=\linewidth]{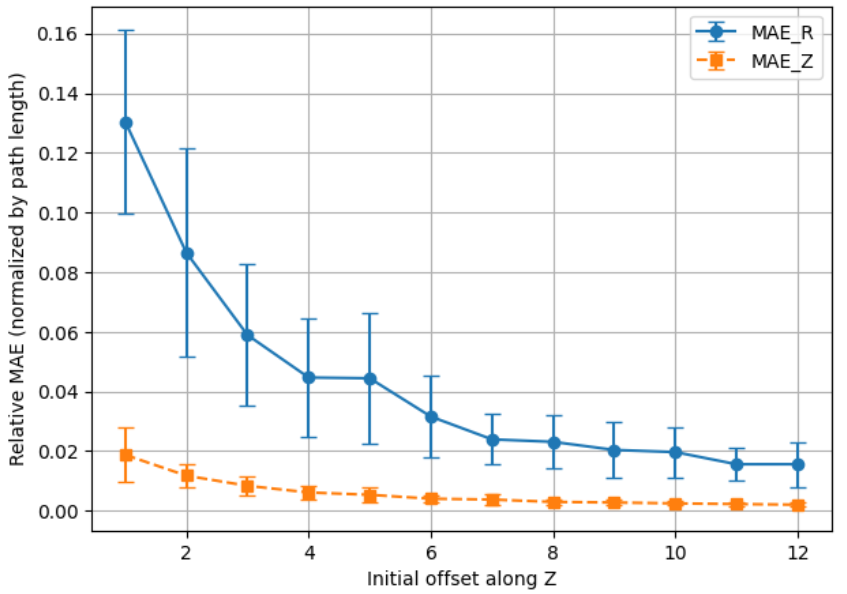}
  \caption{Learning results ($\nu = 1$)}
  \label{Result}
\end{minipage}
\hfill
\begin{minipage}{0.48\textwidth}
  \centering
  \includegraphics[width=\linewidth]{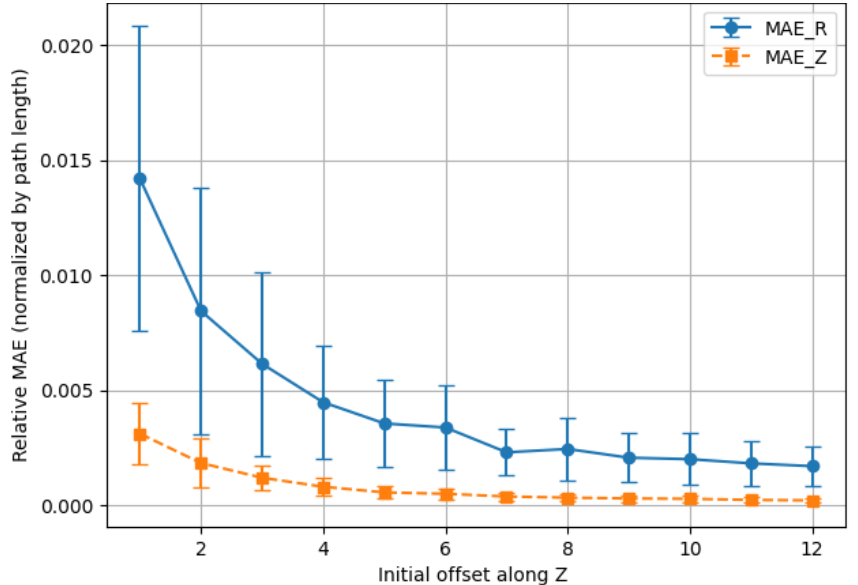}
  \caption{Learning results ($\nu = 0.01$)}
  \label{Result-highReynolds}
\end{minipage}
\end{figure}

\begin{figure}[H]
\begin{center}
\includegraphics[width=0.9\textwidth]{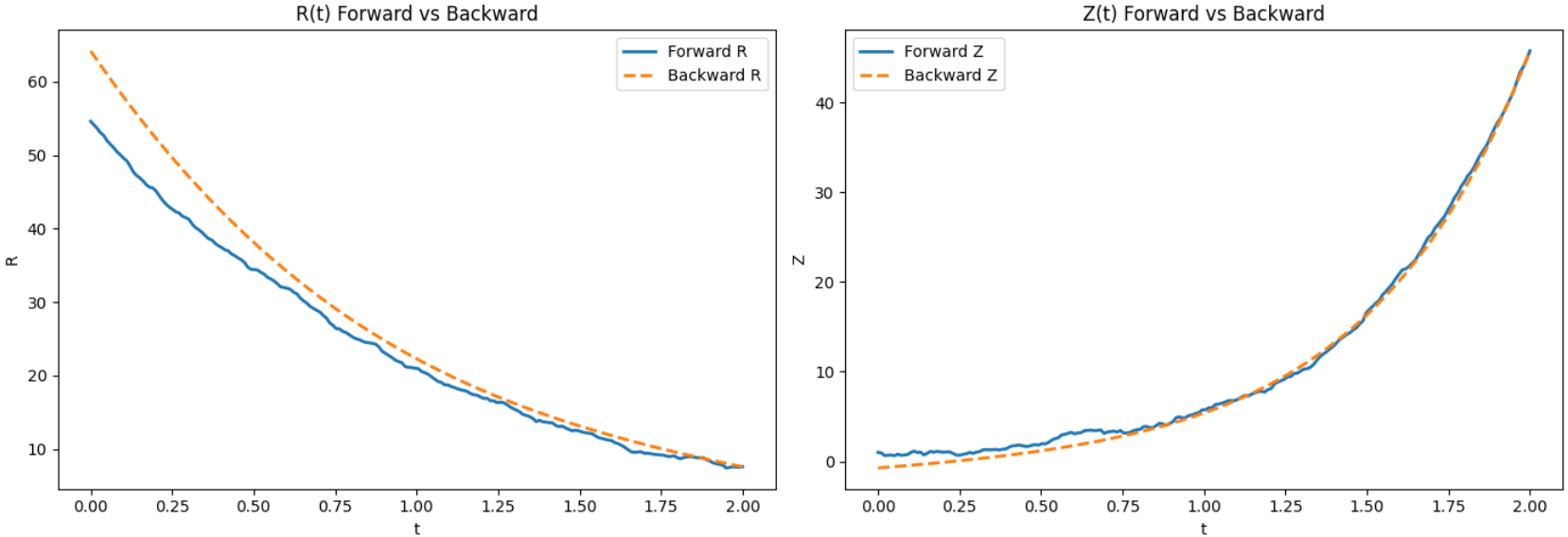}
\caption{Example ($\nu=1, s=1$)}
\end{center}
\end{figure}

\begin{figure}[H]
\begin{center}
\includegraphics[width=0.9\textwidth]{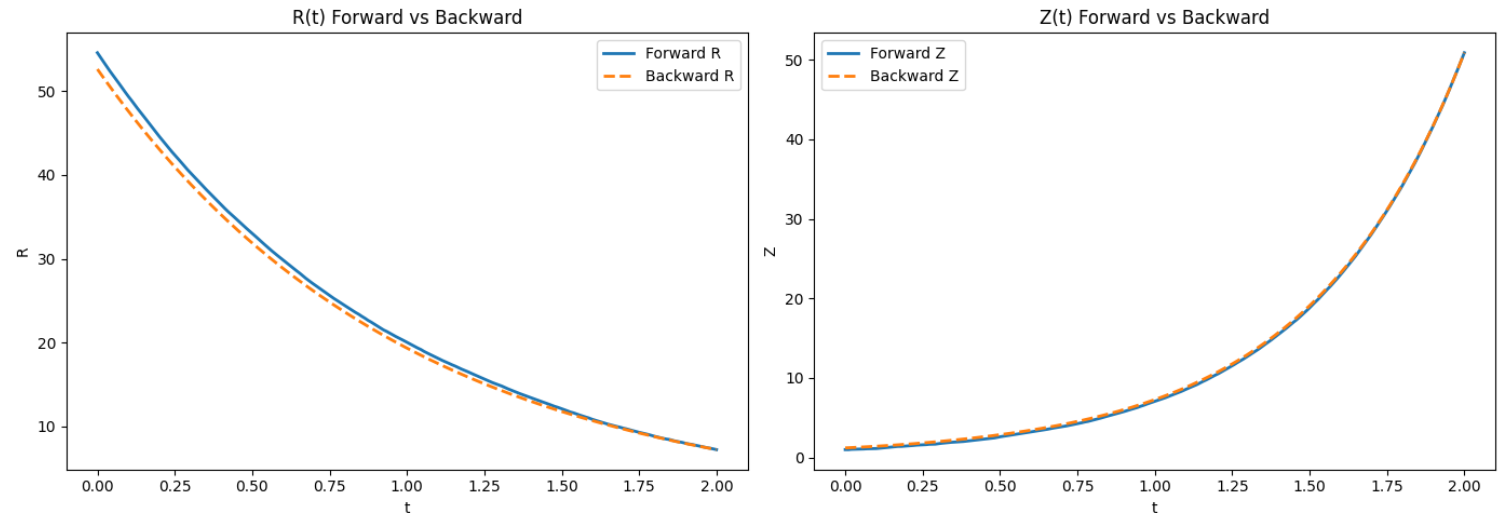}
\caption{Example ($\nu=0.01, s=1$)}
\end{center}
\end{figure}

Comparing Fig.~\ref{Result} ($\nu=1$) and Fig.~\ref{Result-highReynolds} ($\nu=0.01$),
we observe that the graphs exhibit self-similar behavior regardless of the
value of $\nu$.
In both cases, the relative MAE in the $R$ direction is large,
indicating that information about the initial position is significantly
lost as trajectories approach the axis.
In contrast, in the $Z$ direction, the relative MAE remains smaller and
more stable overall (compared with $R$),
indicating that information about the initial position is relatively
well preserved.

\subsection{Additional experiments for the two-dimensional flow}

As an additional experiment, we also carried out numerical simulations
for the following two-dimensional flow ($a=1$)
(using the same notation $R$ and $Z$ for convenience):
\begin{equation*}
\boxed{
\begin{aligned}
R_{k+1} &= -2aR_k\,\Delta t + \sigma\sqrt{\Delta t}
\varepsilon^{(r)},\\[4pt]
Z_{k+1} &= 2aZ_k\,\Delta t + \sigma\sqrt{\Delta t}
\varepsilon^{(z)}.
\end{aligned}
}
\end{equation*}
As in the three-dimensional case, the initial positions of fluid particles
were specified using a scale parameter $s>0$ as
\begin{equation*}
(R_0, Z_0) = (e^{2aT}s, s),
\end{equation*}
and we examined the cases $s = 1,2,\cdots,11,12$.
The learning results are shown below.

\begin{figure}[H]
\centering
\begin{minipage}{0.48\textwidth}
  \centering
  \includegraphics[width=\linewidth]{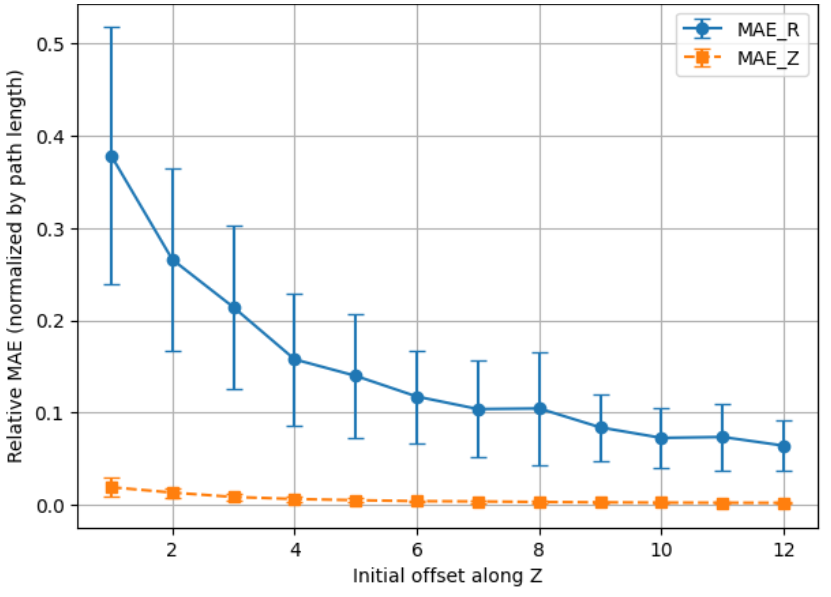}
  \caption{Learning results ($\nu = 1$)}
\end{minipage}
\hfill
\begin{minipage}{0.48\textwidth}
  \centering
  \includegraphics[width=\linewidth]{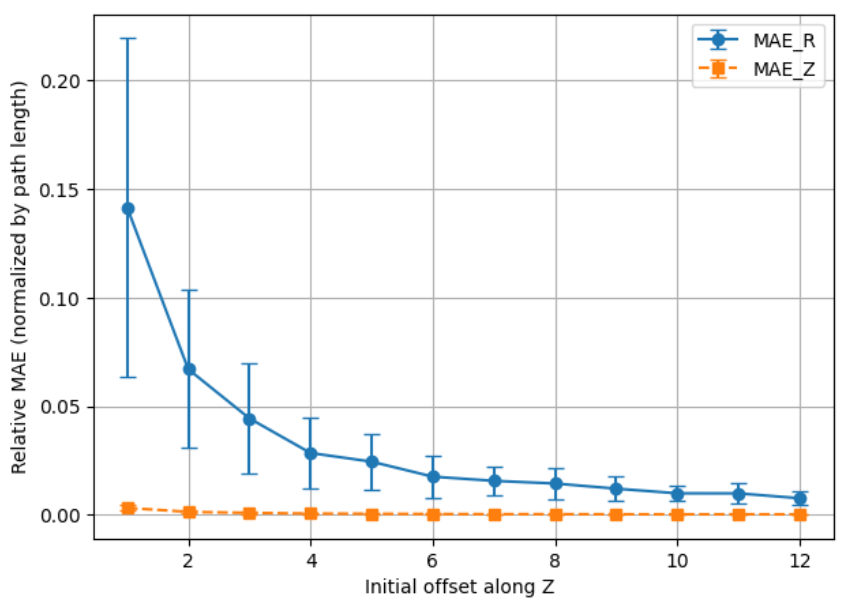}
  \caption{Learning results ($\nu = 0.01$)}
\end{minipage}
\end{figure}

The overall trends are nearly identical to those observed
in the three-dimensional case.

\subsection{Conclusion}

These learning results suggest that

\vspace{0.5cm}

\fbox{
  \begin{minipage}{0.9\linewidth}
In hyperbolic strain fields, the dynamics in the compressive direction
play a role in irreversibly erasing information about the initial position,
whereas such strong information loss is less likely to occur
in the stretching direction.
  \end{minipage}
}

\vspace{0.5cm}

\noindent
As a future direction, it will be important to extend this approach to
more general turbulent flows and to systematically investigate whether
vortex-axis structures and the geometric arrangement of vortices
govern information loss even in settings where axisymmetry is not
explicitly assumed.

\vspace{0.5cm}
\noindent
{\bf Acknowledgments.}\ 
Research of  TY  was partly supported by the JSPS Grants-in-Aid for Scientific
Research 24H00186.

\end{document}